\documentclass[preprint,fleqn,11pt]{elsarticle}
\usepackage{amssymb}
\biboptions{comma,square}

\journal{Journal of Computational and Applied Mathematics}
\usepackage{amssymb,amsmath}
\usepackage{graphics}
\usepackage{graphicx}

\newcommand{\dl}[1]{{\bf Theorem{#1.}}}
\newcommand{\yl}[1]{{\bf Lemma{#1.}}}

\newcommand{\de}[1]{{\bf Definition{#1}}}
\newcommand{\zb}{\textbf{\qquad$\Box$}}
\newcommand{\la}[1]{\label{#1}}
\newcommand{\rf}[1]{(\ref{#1})}

\newcommand{\bo}{\boldsymbol}
\newcommand{\zm}{{\bf Proof.}}

\begin{document}
\thispagestyle{empty} \setcounter{page}{1}

\begin{frontmatter}

\title{A Stable Weak Galerkin Finite
Element Method for Stokes Problem}

\author{Tie Zhang\footnote{Corresponding author at: Department of
Mathematics, Northeastern University, Shenyang 110004, China. {\em
E-mail address} : zhangt@mail.neu.edu.cn(T. Zhang). Tel \& Fax:
+86-024-83680949. },\quad Tao Lin$^\dag$}

\address{$^*$Department of Mathematics, Northeastern University, Shenyang 110004, China\\
$^\dag$Department of Mathematics, Virginia Tech University, Virginia, USA}

\begin{abstract}
We study the weak Galerkin finite element method for
Stokes problem. A new weak Galerkin finite element velocity-pressure space pair is presented which satisfies
the discrete inf-sup condition. Based on this space pair, we establish a stable weak Galerkin approximation scheme without adding any stability term or penalty term. Then, we further derive the optimal error estimates for velocity and pressure approximations, respectively. Numerical experiments are provided to illustrate the theoretical analysis.
\end{abstract}

\end{frontmatter}
\section{Introduction}
\setcounter{section}{1}\setcounter{equation}{0} Recently, the weak
Galerkin finite element method has attracted much attention in the
field of numerical partial differential equations
\cite{Gao,Lin1,Lijian,Lin2,Lin3,Wang,Wang2,Wang3,Wang0,Zhai,Zhang0}. This method was
introduced and analyzed originally in \cite{Wang} for second order elliptic
problems in multi-dimensional domain. In general, a weak Galerkin finite
element method can be considered as an extension of the standard
finite element method where classical derivatives are replaced in
the variational equation by the weakly defined derivatives on
discontinuous weak functions. The main feature of this method is: (1)
the weak derivatives are introduced as distributions for weak functions; (2) the
weak Galerkin finite element function $u_h=\{u_h^0,\,u_h^b\}$ is used in which $u_h^0$ is totally discontinuous on the partition and the value $u_h^b$ of $u_h$ on element edge may be
independent with its value $u_h^0$ in the interior of element. The readers are referred to articles \cite{Lin1,Lin2,Wang2}
for more detailed explanation of this method and its relation with
other finite element methods.

In this paper, we study the weak Galerkin finite element method for Stokes problem.
In the conventional finite element methods solving Stokes and Navier-Stokes
problems, usually the inf-sup condition is required for the velocity-pressure
space pair. The
importance of ensuring the inf-sup condition is widely
understood. Numerical experiments show that the violation of the
inf-sup condition often leads to nonphysical oscillations of the discrete solutions. From
the computational viewpoint, the simple lower-order polynomial space pairs
(for example, the $P_1-P_0$, $Q_1-Q_0$, $P_1-P_1$ and $Q_1-Q_1$
pairs) should be preferred in applications. But unfortunately,
these space pairs do not satisfy the inf$-$sup condition. In
order to circumvent the inf-sup condition, many stabilized methods
were proposed, for example, the penalty methods
and the consistently stabilized
methods, see \cite{Barth,Behr,Girault1,Hughes,Shen}.

For the weak Galerkin finite element methods solving Stokes problem, to the authors' best knowledge, there are only a few articles are presented in existing literatures \cite{Xie,Lin3,Wang3,Zhai}, and all these known methods have a stabilizing term with the penalty factor $h^{-1}$ in the weak Galerkin finite element schemes. However, the stabilizing term or penalty term will add the computation cost and the penalty factor will reduce the numerical stability for $h$ small in solving the discrete linear system. The reasons of adding stabilizing term in known methods are that firstly,  the energy norm $\|\nabla_wu_h\|_h$ of the weak gradient is not a norm on the weak Galerkin spaces used in \cite{Xie,Lin3,Wang3,Zhai}, so the stabilizing terms were added in the weak Galerkin schemes to assure the unique existence of the weak Galerkin solutions; Secondly, the stabilizing term can control the error on element edges so that one can derive the optimal error estimates more easily. In our method here, we choose a different weak Galerkin space from those in \cite{Xie,Lin3,Wang3,Zhai}, for this choice, $\|\nabla_wu_h\|_h$ determines a norm on this weak Galerkin space; furthermore, we use a technique argument to derive the optimal error estimates.

In this paper, we present a weak Galerkin finite element velocity-pressure space pair $V_h\times M_h$ which satisfies
the following discrete inf-sup condition
$$
\sup_{\bo{v}\in V_h}\frac{(\hbox{div}_w\bo{v},q_h)_h}{\|\nabla_w\bo{v}\|_h}\geq \beta\|q_h\|,\;\forall\, q_h\in M_h,
$$
where $\hbox{div}_w$ and $\nabla_w$ are the weak divergence and weak gradient, respectively, see Section 2 for details.
Then, we establish a stable weak Galerkin finite element approximation scheme for the Stokes problem without adding any stability term or penalty term. We adopt a different analysis approach from those in \cite{Xie,Lin3,Wang3,Zhai}.  Using the inf-sup condition and a discrete embedding inequality for weak Galerkin finite element function, we first give the stability estimate for the discrete velocity and pressure. Then, by means of two projection approximations for the velocity and pressure functions, respectively,
we derive the optimal error estimates for the velocity and pressure approximations in various norms.
 We emphasize that our method here also can be applied to solve the Navier-Stokes problem if we deal with the nonlinear convection term properly. This is our ongoing work.

This paper is organized as follows. In Section 2, we first introduce the concepts of weak function, weak gradient and weak divergence, and then we establish the
weak Galerkin finite element approximation for the Stokes problem. Section 3 is devoted to the low order weak Galerkin finite element analysis. Based on some special projection approximations, the stability estimate is established and the optimal error estimates are derived for velocity and pressure approximations, respectively, in the $H^1$ and $L_2$ norms. In Section 4, the method and result in Section 3 are expanded to the high order weak Galerkin finite element approximations. In Section 5, some numerical experiments are provided to illustrate the theoretical analysis. Finally, conclusion is given in Section 6.

Throughout this paper, for an integer $m$, we adopt the notations
$H^{m}(D)$ to indicate the usual Sobolev spaces on domain
$D\subset \Omega$ equipped with the norm $\|\cdot\|_{m,D}$ and
semi-norm $|\cdot|_{m,D}$, and we omit the index $D$ if $D=\Omega$.
The notations $(\cdot,\cdot)$
and $\|\cdot\|$ denote the inner product and norm in the space
$L_2(\Omega)$, respectively. We will use the letter $C$ to represent a
generic positive constant, independent of the mesh size $h$.

\section{Problem and its weak Galerkin finite element approximation}
\setcounter{section}{2}\setcounter{equation}{0}
Consider the Stokes equations
\begin{eqnarray}
-\triangle\bo{u}+\nabla p&=&\bo{
f},\;\;\;in\;\;\Omega,\label{2.1}\\
\hbox{div}\,\bo{u}&=&0,\;\;\;in\;\;\Omega,\label{2.2}\\
\bo{u}&=&\bo{0},\;\;\;on\;\;\partial\Omega,\label{2.3}
\end{eqnarray}
where $\Omega\subset R^2$ is a convex polygonal domain with
boundary $\partial\Omega$, symbols $\triangle$, $\nabla$ and $div$
denote the Laplacian, gradient and divergence operators,
respectively, and $\bo{u}=(u_1,u_2)^T$ represents the velocity, $p$
the pressure and $\bo{f}\in [L_2(\Omega)]^2$ the external volumetric force acting on
the fluid.

The weak form in the primary velocity-pressure formulation for problem \rf{2.1}$\sim$\rf{2.3} is that find $(\bo{u},p)\in [H^1_0(\Omega)]^2\times L^2_0(\Omega)$ such that
\begin{eqnarray}
&&(\nabla\bo{u},\nabla\bo{v})-(p,\hbox{div}\bo{v})=(\bo{f},\bo{v}),\;\forall\,\bo{v}\in [H^1_0(\Omega)]^2,\la{2.4}\\
&&(\hbox{div}\bo{u},q)=0,\;\forall\, q\in L_0^2(\Omega),\la{2.5}
\end{eqnarray}
where the space $L_0^2(\Omega)=\{ q\in L_2(\Omega): \int_{\Omega}qdx=0\}$.

We first introduce the concepts of weak function, weak gradient and weak divergence, and then establish a stable weak Galerkin finite element scheme for Stokes problem \rf{2.1}$\sim$\rf{2.3}.
\subsection{Discrete weak function, weak gradient and weak divergence}
Let $T_h=\bigcup\lbrace K \rbrace$ be a regular triangulation of
domain $\Omega$ so that $\overline{\Omega}=\bigcup_{K\in
T_h}\{K\}$, where the mesh size $h=max \,h_K$, $h_K$ is the
diameter of element $K$.

A weak function on element
$K$ refers to a function $v=\{ v^0,v^b\}$ with $v^0=v|_{K}\in
L_2(K)$ and $v^b=v|_{\partial K}\in L_2(\partial K)$. Note that for
a weak function $v=\{v^0,v^b\}$, $v^b$ may not be necessarily the
trace of $v^0$ on element boundary $\partial K$.

There are many kinds of weak Galerkin finite element spaces that can be used in the numerical PDEs \cite{Wang}. For problem \rf{2.1}$\sim$\rf{2.3}, we introduce the following local weak Galerkin finite element space. For any non-negative integer $l\geq 0$, let $P_l(K)$ be the space composed of all polynomials on $K$ with degree no more than
$l$, and the polynomial space $P_l(\partial K)=\{p: p|_e\in P_l(e), \;\hbox{edge}\,\; e\subset\partial K\}$. Define the weak Galerkin finite element space on element $K$ by
\begin{equation}
W_{k,k+1}(K)=\{\,v=\{v^0,v^b\}: v^0\in P_k(K),\,v^b\in P_{k+1}(\partial
K),\,k\geq 0\}.\la{2.6}
\end{equation}
In analogy with the Green formula held for smooth function, we define the discrete weak gradient and weak divergence of a weak Galerkin finite element function $v\in W_{k,k+1}(K)$ as follows, see \cite{Wang}.\\
\de{ 2.1.}\quad For weak function $v=\{v^0,v^b\}\in W_{k,k+1}(K)$, its discrete weak gradient
$\nabla_wv\in [P_{k+1}(K)]^2$ is defined as the unique solution
of equation:
\begin{equation}
\int_{K}\nabla_{w}v\cdot\boldsymbol{q}dx=-\int_Kv^0\hbox{div}\boldsymbol{q}dx+\int_{\partial
K}v^b\boldsymbol{q}\cdot \bo{n}ds,\,\forall\, \boldsymbol{q}\in
[P_{k+1}(K)]^2,\la{2.7}
\end{equation}
where $\bo{n}=(n_1,n_2)^T$ is the outward unit normal vector on $\partial K$.\\
\de{ 2.2.}\quad For vector weak function $\bo{v}=(v_1,v_2)^T\in [W_{k,k+1}(K)]^2$, its discrete weak divergence
$\hbox{div}_w \bo{v}\in P_{k+1}(K)$ is defined as the unique solution of equation:
\begin{equation}
\int_{K}\hbox{div}_w\bo{v}qdx=-\int_K\bo{v}^0\cdot\nabla qdx+\int_{\partial
K}\bo{v}^b\cdot \bo{n} qds,\,\forall\, q\in
P_{k+1}(K),\la{2.8}
\end{equation}
where $\bo{v}^0=(v^0_1,v^0_2)^T,\,\bo{v}^b=(v^b_1,v^b_2)^T$.

For vector function $\bo{u}=(u_1,u_2)^T$ and $\bo{v}=(v_1,v_2)^T$, and matrix function $\bo{\tau}=(\bo{\tau_1},\bo{\tau}_2)^T$ and $\bo{\sigma}=(\bo{\sigma}_1,\bo{\sigma}_2)^T$, as usual, we set the inner-product and differential operation rule:
\begin{eqnarray*}
&&(\bo{u},\bo{v})=(u_1,v_1)+(u_2,v_2),\;\;(\bo{\tau},\bo{\sigma})=(\bo{\tau}_1,\bo{\sigma}_1)+(\bo{\tau}_2,\bo{\sigma}_2),\\
&&\nabla\bo{u}=(\nabla u_1,\nabla u_2)^T,\;\;\hbox{div}\bo{\tau}=(\hbox{div}\bo{\tau}_1,\hbox{div}\bo{\tau}_2)^T.
\end{eqnarray*}
Then, according to definition \rf{2.7}, for a vector weak function $\bo{v}=(v_1,v_2)^T\in [W_{k,k+1}(K)]^2$, its weak gradient $\nabla_w\bo{v}=(\nabla_wv_1,\nabla_wv_2)^T\in [P_{k+1}(K)]^{2\times 2}$ is the unique solution of equation:
\begin{equation}
\int_{K}\nabla_{w}\bo{v}\cdot\boldsymbol{\tau}dx=-\int_K\bo{v}^0\cdot\hbox{div}\boldsymbol{\tau}dx+\int_{\partial
K}\bo{v}^b\cdot (\boldsymbol{\tau}\bo{n})ds,\,\forall\, \boldsymbol{\tau}\in [P_{k+1}(K)]^{2\times 2}.\la{2.9}
\end{equation}

Two important properties on weak gradient and weak divergence can be stated as follows.\\
\yl{ 2.1}(see \cite[Lemma 5.1]{Wang})\quad {\em Let $v=\{ v^0,v^b\}\in W_{k,k+1}(K)$ be a weak
function. Then, $\nabla_{w}v=0$ on $K$ if and only if
$v=constant$, that is, $v^0=v^b=constant$ on $K$.}\\
\yl{ 2.2}\quad{\em For vector weak function $\bo{v}\in [W_{k,k+1}(K)]^2$, we have}\
\begin{eqnarray}
\|\hbox{div}_w\bo{v}\|_{0,K}\leq \sqrt{2}\,\|\nabla_w\bo{v}\|_{0,K}.\la{2.10}
\end{eqnarray}
\zm\quad For $q\in P_{k+1}(K)$, let diagonal matrix $\bo{\tau}=\hbox{diag}(q,q)$. Then, from definitions \rf{2.9} and \rf{2.8}, we have
\begin{eqnarray*}
\int_{K}\nabla_{w}\bo{v}\cdot\boldsymbol{\tau}dx=-\int_K\bo{v}^0\cdot(\partial_{x_1}q,\partial_{x_2}q)^Tdx+\int_{\partial
K}\bo{v}^b\cdot (qn_1,qn_2)^Tds\\
=-\int_K\bo{v}^0\cdot \nabla qdx+\int_{\partial
K}\bo{v}^b\cdot \bo{n}q ds=\int_{K}\hbox{div}_w\bo{v}qdx,\,\forall\, \boldsymbol{\tau}=\hbox{diag}(q,q).
\end{eqnarray*}
Hence, taking $q=\hbox{div}_w\bo{v}$, we obtain
$$
\|\hbox{div}_w\bo{v}\|_{0,K}^2\leq \|\nabla_w\bo{v}\|_{0,K}\|\bo{\tau}\|_{0,K}\leq \sqrt{2}\,\|\nabla_w\bo{v}\|_{0,K}\|q\|_{0,K}.
$$
The proof is completed.\zb

The following trace inequality will be used frequently in our analysis.
\begin{eqnarray*}
&&\|u\|_{L_2(\partial K)}\leq
Ch_K^{-\frac{1}{2}}\big(\,\|u\|_{0,K}+h_K\|\nabla
u\|_{0,K}\big),\;u\in H^1(K).
\end{eqnarray*}

\subsection{Weak Galerkin finite element scheme}
Introduce the weak Galerkin finite element spaces on triangulation $T_h$:
\begin{eqnarray*}
&&S_h=\{v=\{v^0,v^b\}:\;v|_K\in W_{k,k+1}(K),\,v^b|_{\partial K}\;\hbox{is single valued},\; K\in
T_h\},\\
&&S_h^0=\{v=\{v^0,v^b\}\in S_h:\, v^b|_{\partial\Omega}=0\}.
\end{eqnarray*}
It should be pointed out that for $v=\{v^0,v^b\}\in S_h$, the single value condition of $v^b$ on $\partial K$ implies that $v^b$ is continuous across $\partial K$. On the other hand, the component $v^0$ is defined element-wise and completely discontinuous on $T_h$. A weak Galerkin finite element function $v=\{v^0,v^b\}\in S_h$ is glued in different elements by $v^b$.

Usually, in the finite element analysis, for a properly smooth function $u$, one can find an approximation function in $S_h$ which can approximate $u$ well. But the usual projection function and interpolation function of $u$ are not in the weak function space $S_h$. So we need to introduce a new projection function in $S_h$. For $l\geq 0$, let $P_h^l$ be the local $L_2$ projection operator, restricted on each element $K$, $P_h^l:\,u\in L_2(K)\rightarrow
P_h^lu\in P_l(K)$ such that
\begin{equation}
(u-P_h^lu,q)_{K}=0,\;\forall\,q\in P_l(K),\,K\in T_h.\la{2.11}
\end{equation}
By the Bramble-Hilbert lemma, it is easy to prove that (see
\cite{Zhang})
\begin{equation}
\|u-P_h^lu\|_{0,K}\leq Ch_K^{s}\|u\|_{s,K},\;0\leq s\leq
l+1.\la{2.12}
\end{equation}
Furthermore, let $P_{\partial K}^{k+1}: L_2(e)\rightarrow P_{k+1}(e),\,e\subset \partial K$, be the $L_2$ projection operator. Now,  we define a projection operator $Q_h: u\in H^1(\Omega)\rightarrow Q_hu\in S_h$, restricted on each element $K$,
\begin{equation}
Q_hu|_K=\{Q_h^0u,Q_h^bu\}=\{P^k_hu,P^{k+1}_{\partial K}u\},\;K\in T_h.\la{2.13}
\end{equation}
For a vector function $\bo{u}$, we set $P_h^k\bo{u}=(P_h^ku_1,P_h^ku_2)^T$, $Q_h\bo{u}=(Q_hu_1,Q_hu_2)^T$.

For function $v$ or $\bo{v}$ defined on $T_h$, we set the global operation
$$
(\nabla_wv)|_K=\nabla_w(v|_K),\;(\hbox{div}_w\bo{v})|_K=\hbox{div}_w(\bo{v}|_K),\;K\in T_h.
$$
Projection $Q_hu$ has the following important properties.\\
\yl{ 2.3}\quad{\em Let $\bo{u}\in [H^{1+s}(\Omega)]^2, s\geq 0$. Then, we have}
\begin{eqnarray}
&&\hbox{div}_wQ_h\bo{u}=P_h^{k+1}(\hbox{div}\bo{u}),\;\;\nabla_wQ_h\bo{u}=P_h^{k+1}(\nabla\bo{u}),\la{2.14}\\
&&\|\bo{u}-Q_h^0\bo{u}\|_{0,K}\leq Ch_K^{s}\|\bo{u}\|_{s,K},\,0\leq s\leq
k+1,\;K\in T_h,\la{2.15}\\
&&\|\nabla_{w}Q_h\bo{u}-\nabla\bo{u}\|_{0,K}\leq
Ch^s_K\|\bo{u}\|_{1+s,K},\,0\leq s\leq k+2,\;K\in T_h.\la{2.16}
\end{eqnarray}
\zm\quad From definitions \rf{2.8}, \rf{2.13} and the Green's formula, we have
\begin{eqnarray*}
&&\int_K\hbox{div}_wQ_h\bo{u}qdx=-\int_KQ_h^0\bo{u}\cdot\nabla qdx+\int_{\partial K}Q_h^b\bo{u}\cdot\bo{n} qds\\
&=&-\int_K\bo{u}\cdot\nabla qdx+\int_{\partial K}\bo{u}\cdot\bo{n}qds
=\int_K\hbox{div}\bo{u}qdx,\;\forall\,q\in P_{k+1}(K).
\end{eqnarray*}
This shows that $\hbox{div}_wQ_h\bo{u}=P_h^{k+1}(\hbox{div}\bo{u})$ holds. Similarly, from \rf{2.9} and \rf{2.13}, we can derive  $\nabla_wQ_h\bo{u}=P_h^{k+1}(\nabla\bo{u})$, which also implies approximation property \rf{2.16}. Estimate \rf{2.15} comes from the fact that $Q_h^0=P_h^k$.\zb

Denote the discrete $L_2$ inner product and norm by
$$
(u,v)_h=\sum_{K\in T_h}(u,v)_{K}=\sum_{K\in
T_h}\int_Ku\,vdx,\;\;\;\;\|u\|_h^2=(u,u)_h.
$$
Now, we introduce the velocity and pressure approximation spaces:
$$
V_h=[S^0_h]^2,\;\;M_h=\{q_h\in L^2_0(\Omega): q_h|_K\in P_k(K),\,K\in T_h\}.
$$
Obviously, $Q_hu\in S^0_h$ if $u\in H^1_0(\Omega)$ so that $Q_h\bo{u}\in V_h$ if $\bo{u}\in [H^1_0(\Omega)]^2$.

Motivated by weak form \rf{2.4}--\rf{2.5}, we define the weak Galerkin finite
element approximation of problem \rf{2.1}$\sim$\rf{2.3} by finding $(\bo{u}_h,p_h)\in V_h\times M_h$ such that
\begin{eqnarray}
&&(\nabla_w\bo{u}_h,\nabla_w\bo{v})_h-(p_h,\hbox{div}_w\bo{v})_h=(\bo{f,v}^0),\,\forall\,\bo{v}\in V_h,\la{2.17}\\
&&(\hbox{div}_w\bo{u}_h,q_h)_h=0,\;\forall\,q_h\in M_h.\la{2.18}
\end{eqnarray}
It is well known that the inf-sup condition is very important for the Stokes problem and its finite element analysis. In analogy with the  conventional finite element method, we here also establish a discrete inf-sup condition for the weak Galerkin finite element space pair $V_h\times M_h$.\\
\yl{ 2.4}\quad{\em For weak Galerkin finite element space pair $V_h\times M_h$, the following discrete inf-sup condition holds}
\begin{equation}
\sup_{\bo{v}\in V_h}\frac{(\hbox{div}_w\bo{v},q_h)_h}{\|\nabla_w\bo{v}\|_h}\geq \beta\|q_h\|,\;\forall\, q_h\in M_h,\la{2.19}
\end{equation}
{\em where $\beta$ is a constant independent of $h$}.\\
\zm\quad For any given $q_h\in M_h\subset L^2_0(\Omega)$, it is well known that there exists a function $\bo{w}\in [H^1_0(\Omega)]^2$ and constant $C_0$ such that (see \cite{Girault1})
\begin{equation}
\hbox{div}\bo{w}=q_h,\;\;\|\bo{w}\|_1\leq C_0\|q_h\|.\la{2.20}
\end{equation}
Hence, from Lemma 2.3, we first obtain
$$
\|\nabla_wQ_h\bo{w}\|_h=\|P_h^{k+1}\nabla\bo{w}\|_h\leq \|\nabla\bo{w}\|\leq C_0\|q_h\|,
$$
and then
$$
\frac{(\hbox{div}_wQ_h\bo{w},q_h)_h}{\|\nabla_wQ_h\bo{w}\|_h}=\frac{(P_h^{k+1}(\hbox{div}\bo{w}),q_h)_h}{\|\nabla_wQ_h\bo{w}\|_h}
=\frac{(q_h,q_h)_h}{\|\nabla_wQ_h\bo{w}\|_h}\geq C_0^{-1}\|q_h\|.
$$
This implies inf-sup condition \rf{2.19}.\zb

Introduce the norm notation
\begin{equation}
\|v\|^2_{1,h}=\|\nabla v^0\|^2_h+\sum_{K\in T_h}\int_{\partial K}h_K^{-1}(v^0-v^b)^2ds,\;v\in S^0_h.\la{2.21}
\end{equation}
\yl{ 2.5}\quad{\em Both $\|\nabla_wv\|_h$ and $\|v\|_{1,h}$ are norm on space $S_h^0$ and this two norms are equivalent, that is, there exist positive constants $C_1$ and $C_2$ independent of $h$ such that}
\begin{equation}
  C_1\|v\|_{1,h}\leq \|\nabla_wv\|_h\leq C_2\|v\|_{1,h},\;\forall\,v\in S_h^0.\la{2.22}
\end{equation}
\zm\quad We only need to prove that $\|\nabla_wv\|_h$ is a norm on $S_h^0$ and \rf{2.22} holds. Let $v\in S_h^0$ and $\|\nabla_{w}v\|_h=0$. Then, from Lemma
2.1, we know that $v=\{v^0,v^b\}$ is piecewise constant on $T_h$, that is, $v^0=v^b=constant$ on each element $K$. Since $v^b$ is continuous acrose $\partial K$ and $v^b|_{\partial\Omega}=0$, so we have $v=0\, (v^0=v^b=0)$ which implies $\|\nabla_wv\|_h$ is a norm on $S_h^0$. The equivalence demonstration of norms $\|\nabla_wv\|_h$ and $\|v\|_{1,h}$ can be found in \cite[Lemma 3.2]{Mu}.\zb

By means of Lemma 2.4 and Lemma 2.5, we can obtain the following result. \\
\dl{ 2.1}\quad{\em Weak finite element equations \rf{2.17}-\rf{2.18} has a unique solution $(\bo{u}_h,$ $ p_h) \in V_h\times M_h$}.\\
\zm\quad Since equations \rf{2.17}-\rf{2.18} is a linear system
of equations, we only need to prove the uniqueness. Let $\bo{f}=\bo{0}$,
we need to prove $\bo{u}_h=p_h=0$. Taking $\bo{v=u}_h$ in \rf{2.17} and using \rf{2.18}, we
obtain $\|\nabla_{w}\bo{u}_h\|_h=0$ which together with Lemma 2.5 imply that $\bo{u}_h=\bo{0}$ on $T_h$. Furthermore, with $\bo{u}_h=\bo{f}=\bo{0}$ in equation \rf{2.17}, we have
$$
(p_h,\hbox{div}_w\bo{v})_h=0,\,\forall\,\bo{v}\in V_h.
$$
The proof is completed by using the inf-sup condition \rf{2.19}. \zb

\section{Stability and error analysis}
\setcounter{section}{3}\setcounter{equation}{0}
In order to highlight our analysis method and simplify the argument, in this section, we only discuss the low-order ($k=0$) weak Galerkin finite element scheme \rf{2.17}-\rf{2.18}. The high-order method ($k>0$) will be discussed in next section.

In $k=0$ case, the corresponding spaces are as follows.
\begin{eqnarray*}
&&S^0_h=\{v:\;v|_K\in W_{0,1}(K),\,v^b|_{\partial K}\;\hbox{is single valued},\; K\in
T_h,\, v^b|_{\partial\Omega}=0\},\\
&&V_h=[S^0_h]^2,\;\;M_h=\{q_h\in L_2^0(\Omega): q_h|_K\in P_0(K),\,K\in T_h\}.
\end{eqnarray*}
Moreover, please bear in mind that all results in Section 2 maintain to hold for $k=0$, for example (see Lemma 2.3),
$$
\hbox{div}_wQ_h\bo{u}=P_h^1(\hbox{div}\bo{u}),\;\;\nabla_wQ_h\bo{u}=P_h^1(\nabla\bo{u}),\;\bo{u}\in [H^1(\Omega)]^2.
$$

\subsection{Stability estimate}
We have proved that the weak Galerkin finite element scheme is stable, that is, the solution $(\bo{u}_h,p_h)$ uniquely exists, but we do not give a stability estimate for solution $(\bo{u}_h,p_h)$. In this subsection, we do the stability estimate.

We first introduce a special projection function $\pi_h\bo{u}$, see \cite{Brezzi}.

Let $e_i$ ($i=1,2,3$) be the edge of element $K$ and space
$H(\hbox{div};\Omega)=\{\boldsymbol{u}\in [L_2(\Omega)]^2:
\hbox{div}\boldsymbol{u}\in L_2(\Omega)\}$.  For function $\phi$, the $curl$ operator is defined by curl\,$\phi=(\partial_{x_2}\phi,-\partial_{x_1}\phi)^T$.

Define the projection operator $\pi_h:
H(\hbox{div};\Omega)\rightarrow H(\hbox{div};\Omega)$, restricted
on element $K\in T_h$, $\pi_h\boldsymbol{u}\in [P_{1}(K)]^2$ satisfies
\begin{eqnarray}
\int_{e_i}(\boldsymbol{u}-\pi_h\boldsymbol{u})\cdot\bo{n}
qds=0,\;\forall\,q\in P_{1}(e_i),\,i=1,2,3.\la{3.1}
\end{eqnarray}
\yl{ 3.1}\quad{\em For $\boldsymbol{u}\in H(\hbox{div}\,;\Omega)$,
the projection $\pi_h\boldsymbol{u}$ uniquely exists and
satisfies}
\begin{eqnarray}
&&(\hbox{div}(\boldsymbol{u}-\pi_h\boldsymbol{u}),q)_K=0,\;\forall\,q\in
P_{0}(K),\;K\in T_h.\la{3.2}
\end{eqnarray}
Furthermore, if $\boldsymbol{u}\in [H^{1+s}(\Omega)]^2,\,0\leq s\leq 1$, then we have
\begin{eqnarray}
&&\|\pi_h\boldsymbol{u}\|_{0,K}\leq
C\|\boldsymbol{u}\|_{1,K},\;K\in
T_h,\la{3.3}\\
&&\|\boldsymbol{u}-\pi_h\boldsymbol{u}\|_{0,K}\leq
Ch_K^{1+s}\|\boldsymbol{u}\|_{1+s,K},\;0\leq s\leq 1,\;K\in
T_h.\la{3.4}
\end{eqnarray}
\zm\quad We first prove the unique existence of
$\pi_h\boldsymbol{u}$. For given $\bo{u}$, the six equations in \rf{3.1} form a
consistent linear system of equations on unknown vector polynomial $\pi_h\bo{u}$, so we only need to prove
that $\pi_h\boldsymbol{u}=0$ if $\boldsymbol{u}=0$. Let $\boldsymbol{u}=0$ in \rf{3.1}. Then, we have
$\pi_h\boldsymbol{u}\cdot\bo{n}=0$ on $\partial K$ and
\begin{eqnarray*}
(\hbox{div}\pi_h\boldsymbol{u},q)_K=\int_{\partial K}\pi_h\boldsymbol{u}\cdot
\bo{n}qds=0,\;\forall\,q\in P_{0}(K).
\end{eqnarray*}
Hence, we obtain $\hbox{div}\pi_h\boldsymbol{u}=0$ on $K$ and $\pi_h\boldsymbol{u}\cdot\bo{n}=0$ on $\partial K$. So there
exists a function $\phi\in P_{2}(K)$ so that
curl$\,\phi=\pi_h\boldsymbol{u}$ (see \cite{Girault}). Since the
tangential derivative $\partial_\tau\phi=\hbox{curl}\,\phi\cdot\bo{n}=\pi_h\boldsymbol{u}\cdot
\bo{n}=0$ on $\partial K$, so $\phi=\phi_0=constant$ on $\partial K$.
Let $p=\phi-\phi_0$. Since $p\in P_{2}(K)$ and $p|_{\partial K}=0$,
then there must be $p=0$ so that $\pi_h\boldsymbol{u}=\hbox{curl}\,p=0$.

Next, we prove conclusions \rf{3.2}$\sim$\rf{3.4}. Equation \rf{3.2} comes
directly from the $div$-formula and equation \rf{3.1}. Moreover,
from the solution representation of linear system of equations
\rf{3.1} and the trace inequality, it is easy to see that on the reference
element $\hat{K}$,
\begin{equation}
\|\hat{\pi}_h\hat{\boldsymbol{u}}\|_{0,\hat{K}}\leq
\hat{C}\|\hat{\boldsymbol{u}}\|_{0,\partial
\hat{K}}\leq\hat{C}(\|\hat{\boldsymbol{u}}\|_{0,\hat{K}}+\|\hat{\nabla}\hat{\boldsymbol{u}}\|_{0,\hat{K}}).\la{3.5}
\end{equation}
Then, \rf{3.3} follows
from \rf{3.5} and a scale argument between $\hat{K}$ and $K$.
From \rf{3.3} and the unique existence, we also obtain
$$
\pi_h\bo{u}=\bo{u},\,\forall
\,\bo{u}\in [P_1(K)]^2,\;\;\|\pi_h\boldsymbol{u}\|_{0,K}\leq
C\|\boldsymbol{u}\|_{1+s,K},\, 0\leq s\leq 1.
$$
Then, estimate \rf{3.4} follows from the Bramble-Hilbert lemma.\zb

Afterwards, for matrix function $\bo{\tau}=(\bo{\tau}_1,\bo{\tau}_2)^T$, we set $\pi_h\bo{\tau}=(\pi_h\bo{\tau}_1,\pi_h\bo{\tau}_2)^T$.

The following discrete embedding inequality is an analogy of the
Poincar\'e inequality in $H^1_0(\Omega)$.\\
\yl{ 3.2}\quad{\em Let $\Omega$ be a polygonal domain. Then, for weak function $v\in S_h^0$, there exists a positive constant $C_1$
independent of $h$ such that}
\begin{equation}
\|v^0\|\leq C_1\|\nabla_{w}v\|_h,\;\forall\, v\in
S_h^0.\la{3.6}
\end{equation}
\zm\quad For $v\in S^0_h$, we first make a smooth domain
$\Omega'\supset \Omega$ ( if $\Omega$ is convex, we may set
$\Omega'=\Omega$) and extend $v^0$ to domain $\Omega'$ by setting
$v^0|_{\Omega'\backslash\Omega}=0$. Then, there exists a function
$\psi\in H^1_0(\Omega')\bigcap H^2(\Omega')$ such that
\begin{eqnarray*}
-\triangle \psi=v^0,\;in\,\,\Omega',\;\;\|\psi\|_{2,\Omega'}\leq
C\|v^0\|.
\end{eqnarray*}
Now we set $\boldsymbol{w}=-\nabla \psi$, then $\boldsymbol{w}\in
[H^1(\Omega)]^2$ satisfies
$$
\hbox{div}\boldsymbol{w}=v^0,\;in\,\,\Omega,\;\;\|\boldsymbol{w}\|_1\leq
\|\psi\|_{2,\Omega'}\leq C\|v^0\|.
$$
Hence, we have from \rf{3.2}, \rf{3.3} and definition \rf{2.7} that
\begin{eqnarray*}
\|v^0\|^2&=&(\hbox{div}\boldsymbol{w},v^0)=(\hbox{div}\pi_h\boldsymbol{w},v^0)\\
&=&\sum_{K\in T_h}\big(-\int_K\nabla_{w}v\cdot
\pi_h\boldsymbol{w}dx+\int_{\partial
K}v^b\pi_h\boldsymbol{w}\cdot\bo{n}ds\big)\\
&=&\sum_{K\in
T_h}-\int_K\nabla_{w}v\cdot\pi_h\boldsymbol{w}dx\leq
\|\nabla_{w}v\|_h\|\pi_h\boldsymbol{w}\|\\
&\leq& C\|\nabla_{w}v\|_h\|\boldsymbol{w}\|_1\leq
C\|\nabla_{w}v\|_h\|v^0\|,
\end{eqnarray*}
where we have used the fact that
\begin{equation}
\sum_{K\in T_h}\int_{\partial K}v^b\pi_h\boldsymbol{w}\cdot
\bo{n}ds=\sum_{K\in T_h}\int_{\partial K}v^b\boldsymbol{w}\cdot
\bo{n}ds=0.\la{3.7}
\end{equation}
The proof is completed.\zb

A direct application of Lemma 3.2 is the following stability estimate.\\
\dl{ 3.1}\quad{\em Let $(\bo{u}_h,p_h)\in V_h\times M_h$ be the solution of weak Galerkin finite element equations
\rf{2.17}-\rf{2.18}. Then we have}
\begin{equation}
\|\nabla_{w}\bo{u}_h\|_h+\|p_h\|\leq
3C_1\|\bo{f}\|.\la{3.8}
\end{equation}
\zm\quad Taking $\bo{v}=\bo{u}_h$ in \rf{2.17} and using \rf{2.18} and \rf{3.6}, we obtain
\begin{equation}
\|\nabla_w\bo{u}_h\|_h^2=(\bo{f,u}_h^0)\leq \|\bo{f}\|\,\|\bo{u}_h^0\|\leq C_1\|\bo{f}\|\,\|\nabla_w\bo{u}_h\|_h.\la{3.9}
\end{equation}
This gives the estimate of $\nabla_w\bo{u}_h$. Furthermore, we have from \rf{2.17}, \rf{3.6} and \rf{3.9} that
\begin{eqnarray*}
(p_h,\hbox{div}_w\bo{v})_h=(\nabla_w\bo{u}_h,\nabla_w\bo{v})_h-(\bo{f,v}^0)\\
\leq 2C_1\|\bo{f}\|\,\|\nabla_w\bo{v}\|_h,\,\forall\,\bo{v}\in V_h.
\end{eqnarray*}
The proof is completed by using the inf-sup condition \rf{2.19}.\zb

\subsection{Error analysis}
In this subsection, we do the error analysis. We first set two bilinear forms
\begin{eqnarray}
&&l_1(\bo{u},\bo{v})=(\nabla_wQ_h\bo{u}-\pi_h(\nabla \bo{u}),\nabla_w\bo{v})_h,\la{3.16}\\
&&l_2(p,\bo{v})=\sum_{K\in T_h}\int_{\partial K}(\bo{v}^0-\bo{v}^b)\cdot\bo{n}(Q_h^0p-p)ds.\la{3.17}
\end{eqnarray}
By means of projections $\pi_h$ and $Q_h$, we can derive the following important equation.\\
\yl{ 3.3}\quad{\em Let $(\bo{u},p)\in [H^2(\Omega)]^2\times H^1(\Omega)$ be the solution of Stokes problem \rf{2.1}-\rf{2.3}. Then, $(Q_h\bo{u},Q_h^0p)\in V_h\times M_h$ satisfies}
\begin{eqnarray}
\left\{\begin{array}{ll}
(\nabla_wQ_h\bo{u},\nabla_w\bo{v})_h-(Q_h^0p,\hbox{div}_w\bo{v})_h=(\bo{f,v}^0)
+l_1(\bo{u},\bo{v})+l_2(p,\bo{v}),\\
(\hbox{div}_wQ_h\bo{u},q_h)_h=0,\;\forall\,(\bo{v},q_h)\in V_h\times M_h.\la{3.18}
\end{array}
\right.
\end{eqnarray}
\zm\quad Let $\bo{v}=\{\bo{v}^0,\bo{v}^b\}\in V_h$ so that on each element $K$, $\bo{v}^0\in [P_0(K)]^2,\,\bo{v}^b\in [P_1(\partial K)]^2$. First, by weak gradient definition \rf{2.9} and properties \rf{3.1}-\rf{3.2} of projection $\pi_h\bo{u}$, we have
\begin{eqnarray}
&&(\nabla_w\bo{v},\pi_h(\nabla\bo{u}))_h=-(\hbox{div}\pi_h(\nabla\bo{u}),\bo{v}^0)_h+\sum_{K\in T_h}\int_{\partial K}\bo{v}^b\cdot(\pi_h(\nabla\bo{u})\bo{n})ds\nonumber\\
&=&-(\hbox{div}(\nabla\bo{u}),\bo{v}^0)_h+\sum_{K\in T_h}\int_{\partial K}\bo{v}^b\cdot(\nabla\bo{u}\bo{n})ds=-(\hbox{div}(\nabla\bo{u}),\bo{v}^0)_h,\la{3.19}
\end{eqnarray}
where we have used the fact that $\bo{v}^b$ is continuous across $\partial K$ and $\bo{v}^b|_{\partial\Omega}=0$. Next, by weak divergence definition \rf{2.8} and the Green's formula, and noting that $\bo{v}^b$ is continuous across $\partial K$ and $\bo{v}^b|_{\partial \Omega}=0$, we have
\begin{eqnarray}
&&-(Q_h^0p,\hbox{div}_w\bo{v})_h=(\bo{v}^0,\nabla Q^0_hp)_h-\sum_{K\in T_h}\int_{\partial K}\bo{v}^b\cdot\bo{n}Q^0_hpds\nonumber\\
&=&(\bo{v}^0,\nabla (Q^0_hp-p))_h+(\bo{v}^0,\nabla p)_h-\sum_{K\in T_h}\int_{\partial K}\bo{v}^b\cdot\bo{n}(Q^0_hp-p)ds\nonumber\\
&=&-(\hbox{div}\bo{v}^0,Q^0_hp-p)_h+\sum_{K\in T_h}\int_{\partial K}(\bo{v}^0-\bo{v}^b)\cdot\bo{n}(Q^0_hp-p)ds+(\bo{v}^0,\nabla p)_h,\nonumber\\
&=&\sum_{K\in T_h}\int_{\partial K}(\bo{v}^0-\bo{v}^b)\cdot\bo{n}(Q^0_hp-p)ds+(\bo{v}^0,\nabla p)_h.\la{3.20}
\end{eqnarray}
Combining \rf{3.19} and \rf{3.20}, and using equation \rf{2.1}, we obtain
$$
(\pi_h(\nabla\bo{u}),\nabla_w\bo{v})_h-(Q^0_hp,\hbox{div}_w\bo{v})_h-l_2(p,\bo{v})=-(\hbox{div}(\nabla\bo{u}),\bo{v}^0)_h+(\nabla p,\bo{v}^0)_h=(\bo{f,v}^0).
$$
Together with \rf{3.16}, this gives the first equation in \rf{3.18}. The second equation in \rf{3.18} comes from the fact that $\hbox{div}Q_h\bo{u}=P_h^1(\hbox{div}\bo{u})=0$.\zb

Now, we can give the optimal error estimates for velocity and pressure approximations.\\
\dl{ 3.2}\quad{\em Let $(\bo{u},p)\in [H^2(\Omega)]^2\times H^1(\Omega)$ and $(\bo{u}_h,p_h)\in V_h\times M_h$ be the solutions of Stokes problem \rf{2.1}-\rf{2.3} and weak Galerkin finite element equations \rf{2.17}-\rf{2.18}, respectively. Then, we have}
\begin{equation}
\|\nabla\bo{u}-\nabla_w\bo{u}_h\|_h+\|p-p_h\|\leq Ch(\|\bo{u}\|_2+\|p\|_1).\la{3.21}
\end{equation}
\zm\quad Let error functions $\bo{e}_h=Q_h\bo{u}-\bo{u}_h,\, \rho_h=Q_h^0p-p_h$. From equations \rf{2.17}-\rf{2.18} and \rf{3.18}, we see that $(\bo{e}_h,\rho_h)\in V_h\times M_h$ satisfies
\begin{eqnarray}
&&(\nabla_w\bo{e}_h,\nabla_w\bo{v})_h-(\rho_h,\hbox{div}_w\bo{v})_h=l_1(\bo{u},\bo{v})+l_2(p,\bo{v}),\forall\,\bo{v}\in V_h,\la{3.22}\\
&&(\hbox{div}_w\bo{e}_h,q_h)_h=0,\;\forall\,q_h\in M_h.\la{3.23}
\end{eqnarray}
From \rf{3.16} and the approximation property of $\pi_h\bo{u}$, we obtain
\begin{eqnarray}
&&|l_1(\bo{u},\bo{v})|\leq (\|\nabla_wQ_h\bo{u}-\nabla\bo{u}\|_h+\|\nabla \bo{u}-\pi_h(\nabla \bo{u})\|_h)\|\nabla_w\bo{v}\|_h\nonumber\\
&\leq& (\|P_h^1(\nabla \bo{u})-\nabla\bo{u}\|_h+Ch\|\nabla \bo{u}\|_1)\|\nabla_w\bo{v}\|_h\leq Ch\|\bo{u}\|_2\|\nabla_w\bo{v}\|_h,\la{3.24}
\end{eqnarray}
Next, using \rf{3.17}, inverse inequality and Lemma 2.5, we have
\begin{eqnarray}
|l_2(p,\bo{v})|\leq C\|\bo{v}\|_{1,h}\|Q_h^0p-p\|\leq Ch\|p\|_1\|\nabla_w\bo{v}\|_h.\la{3.25}
\end{eqnarray}
Substituting estimates \rf{3.24}-\rf{3.25} into \rf{3.22} with $\bo{v}=\bo{e}_h$ and using \rf{3.23}, it yields
\begin{equation}
\|\nabla_w\bo{e}_h\|^2_h\leq Ch(\|\bo{u}\|_2+\|p\|_1)\|\nabla_w\bo{e}_h\|_h.\la{3.25b}
\end{equation}
This gives the velocity estimate by using the triangle inequality and noting that $\nabla_wQ_h\bo{u}=P_h^1(\nabla\bo{u})$. Furthermore, from equation \rf{3.22} and estimates \rf{3.24}-\rf{3.25}, we also obtain
\begin{eqnarray*}
(\rho_h,\hbox{div}_w\bo{v})_h=(\nabla_w\bo{e}_h,\nabla_w\bo{v})_h-l_1(\bo{u},\bo{v})-l_2(p,\bo{v})\\
\leq Ch(\|\bo{u}\|_2+\|p\|_1)\|\nabla_w\bo{v}\|_h,\;\forall\,\bo{v}\in V_h.
\end{eqnarray*}
The pressure estimate is derived by using the inf-sup condition \rf{2.19}.\zb

Below we do error estimate in the $L_2$-norm. To this end, we introduce the auxiliary problem: $(\bo{w},\xi)\in [H^1_0(\Omega)]^2\times L_0^2(\Omega)$ such that \cite{Girault1}
\begin{eqnarray}
&&-\triangle\bo{w}+\nabla \xi=\bo{
e}_h^0,\;\;\;in\;\;\Omega,\;\;\|\bo{w}\|_2+\|\xi\|_1\leq C\|\bo{e}_h^0\|,\label{3.26}\\
&&\hbox{div}\,\bo{w}=0,\;\;\;in\;\;\Omega,\;\;\bo{w}=\bo{0},\;\;\;on\;\;\partial\Omega,\label{3.27}
\end{eqnarray}
where error functions $\bo{e}_h=Q_h\bo{u}-\bo{u}_h$.\\
\dl{ 3.3}\quad{\em Let $(\bo{u},p)\in [H^2(\Omega)]^2\times H^1(\Omega)$ and $(\bo{u}_h,p_h)\in V_h\times M_h$ be the solutions of Stokes problem \rf{2.1}-\rf{2.3} and weak Galerkin finite element equations \rf{2.17}-\rf{2.18}, respectively. Then, we have}
\begin{equation}
\|Q_h^0\bo{u}-\bo{u}^0_h\|\leq Ch^2(\|\bo{u}\|_2+\|p\|_1)+Ch\|\bo{f}-Q_h^0\bo{f}\|.\la{3.28}
\end{equation}
\zm\quad Let $(\bo{w},\xi)$ be the solution of problem \rf{3.26}-\rf{3.27}. By a similar argument to that of Lemma 3.4, we see that $(\bo{w},\xi)$ satisfies (see \rf{3.18})
\begin{equation}
(\nabla_wQ_h\bo{w},\nabla_w\bo{v})_h-(Q_h^0\xi,\hbox{div}_w\bo{v})_h=(\bo{e}_h^0,\bo{v}^0)
+l_1(\bo{w},\bo{v})+l_2(\xi,\bo{v}),\,\bo{v}\in V_h.\nonumber
\end{equation}
Taking $\bo{v}=\bo{e}_h$, it follows from error equation \rf{3.23} that
\begin{equation}
\|\bo{e}_h^0\|^2=(\nabla_wQ_h\bo{w},\nabla_w\bo{e}_h)_h-l_1(\bo{w},\bo{e}_h)-l_2(\xi,\bo{e}_h).\la{3.29}
\end{equation}
From estimates \rf{3.24}-\rf{3.25}, we obtain
\begin{eqnarray}
&&|l_1(\bo{w},\bo{e}_h)|\leq Ch\|\bo{w}\|_2\|\nabla_w\bo{e}_h\|_h,\la{3.30}\\
&&|l_2(\xi,\bo{e}_h)|\leq Ch\|\xi\|_1\|\nabla_w\bo{e}_h\|_h.\la{3.31}
\end{eqnarray}
Below we estimate the first term in \rf{3.29}. By using equations \rf{2.4} and \rf{2.17} satisfied by $\bo{u}$ and $\bo{u}_h$, respectively, and noting that $\hbox{div}_wQ_h\bo{w}=P_h^1(\hbox{div}\bo{w})=0$, we obtain
\begin{eqnarray*}
&&(\nabla_wQ_h\bo{w},\nabla_w\bo{e}_h)_h=(\nabla_wQ_h\bo{w},P_h^1(\nabla\bo{u})-\nabla_w\bo{u}_h)_h\\
&=&(\nabla_wQ_h\bo{w}, \nabla\bo{u})_h-(\nabla_wQ_h\bo{w},\nabla_w\bo{u}_h)_h\\
&=&(P_h^1(\nabla\bo{w})-\nabla\bo{w}, \nabla\bo{u})_h+(\nabla\bo{u},\nabla\bo{w})-(\nabla_w\bo{u}_h,\nabla_wQ_h\bo{w})_h\\
&=&(P_h^1(\nabla\bo{w})-\nabla\bo{w}, \nabla\bo{u}-P_h^1(\nabla\bo{u}))_h+(\bo{f},\bo{w})-(p,\hbox{div}\bo{w})\\
&&-(\bo{f},Q_h^0\bo{w})-(p_h,\hbox{div}_wQ_h\bo{w})_h\\
&=&(P_h^1(\nabla\bo{w})-\nabla\bo{w}, \nabla\bo{u}-P_h^1(\nabla\bo{u}))_h+(\bo{f}-Q_h^0\bo{f},\bo{w}-Q_h^0\bo{w}).
\end{eqnarray*}
Hence, we have from the approximation properties that
\begin{equation}
(\nabla_wQ_h\bo{w},\nabla_w\bo{e}_h)_h\leq C(h^2\|\bo{u}\|_2+h\|\bo{f}-Q_h^0\bo{f}\|)\|\bo{w}\|_2.\la{3.32}
\end{equation}
Substituting estimates \rf{3.30}$\sim$\rf{3.32} into \rf{3.29} and using \rf{3.25b}, it yields
$$
\|\bo{e}^0_h\|^2\leq C\Big(h^2(\|\bo{u}\|_2+\|p\|_1)+h\|\bo{f}-Q_h^0\bo{f}\|\Big)(\|\bo{w}\|_2+\|\xi\|_1).
$$
This gives the desired estimate by using \rf{3.26}.\zb

If $\bo{f}\in [H^1(\Omega)]^2$, from \rf{3.28} we also obtain
\begin{equation}
\|Q_h^0\bo{u}-\bo{u}^0_h\|\leq Ch^2(\|\bo{u}\|_2+\|p\|_1+\|\bo{f}\|_1).\la{3.33}
\end{equation}
The regularity requirements for $\bo{u}, p$ and $\bo{f}$ in \rf{3.33} are the same as that of the finite volume element method for Stokes problem.
\section{High order weak Galerkin finite element approximation}
\setcounter{section}{4}\setcounter{equation}{0} If the solutions of Stokes equations \rf{2.1}$\sim$\rf{2.3} have higher regularity, for example, $(\bo{u},p)\in [H^{2+k}(\Omega)]^2\times H^{1+k}(\Omega),\,k\geq 1$, we may consider to use the high order weak Galerkin finite element spaces $V_h\times M_h$ with $k\geq 1$. In this section, we will prove that the weak Galerkin finite element scheme \rf{2.17}-\rf{2.18} with $k\geq 1$ still works well and the corresponding optimal error estimates maintain to hold.

To extend our analysis and results to high order weak Galerkin finite element scheme, the only task needed to be done is to extend this projection $\pi_h$ to the high order polynomial space, the remanent arguments are completely parallel to those in Section 3.

Let $e_i$ and $\lambda_i$ ($1\leq i\leq 3$) be the edge and
barycenter coordinate of $K$, respectively. Let space
$P_{k+2}^0(K) =\{\,p\in P_{k+2}(K): p|_{\partial
K}=0\}=\lambda_1\lambda_2\lambda_3P_{k-1}(K)$.
For $k\geq 1$, we define the projection operator $\pi_h:
H(\hbox{div};\Omega)\rightarrow H(\hbox{div};\Omega)$, restricted
on each element $K\in T_h$, $\pi_h\boldsymbol{u}\in [P_{k+1}(K)]^2$ satisfies
\begin{eqnarray}
&&(\boldsymbol{u}-\pi_h\boldsymbol{u},\nabla
q)_K=0,\;\forall\,q\in
P_k(K),\la{4.1}\\
&&\int_{e_i}(\boldsymbol{u}-\pi_h\boldsymbol{u})\cdot\bo{n}
qds=0,\;\forall\,q\in P_{k+1}(e_i),\,i=1,2,3,\la{4.2}\\
&&(\boldsymbol{u}-\pi_h\boldsymbol{u},\hbox{
curl}\,q)_K=0,\;\forall\, q\in P^0_{k+2}(K).\la{4.3}
\end{eqnarray}
Some properties of projection $\pi_h\boldsymbol{u}$ had been
discussed in \cite{Brezzi}, we here give a more detailed analysis
for our argument requirement.\\
\yl{ 4.1}\quad{\em For $\boldsymbol{u}\in H(\hbox{div}\,;\Omega)$,
the projection $\pi_h\boldsymbol{u}$ uniquely exists and
satisfies}
\begin{eqnarray}
&&(\hbox{div}(\boldsymbol{u}-\pi_h\boldsymbol{u}),q)_K=0,\;\forall\,q\in
P_{k}(K),\;K\in T_h.\la{4.4}
\end{eqnarray}
Furthermore, if $\boldsymbol{u}\in [H^{1+s}(\Omega)]^2,\,s\geq 0$, then
\begin{eqnarray}
&&\|\pi_h\boldsymbol{u}\|_{0,K}\leq
C\|\boldsymbol{u}\|_{1,K},\;K\in
T_h,\la{4.5}\\
&&\|\boldsymbol{u}-\pi_h\boldsymbol{u}\|_{0,K}\leq
Ch_K^{1+s}\|\boldsymbol{u}\|_{1+s,K},\;0\leq s\leq k+1,\;K\in
T_h.\la{4.6}
\end{eqnarray}
\zm\quad We first prove the unique existence of
$\pi_h\boldsymbol{u}$. Since the number of dimensions (noting that
\rf{4.1} is trivial for $q=constant$):
$$
\hbox{dim}(P_k(K))-1+3\,\hbox{dim}(P_{k+1}(e_i))+\hbox{dim}(P^0_{k+2}(K))=2\,\hbox{dim}(P_{k+1}(K)),
$$
so the linear system of equations \rf{4.1}$\sim$\rf{4.3} is
consistent. Thus, we only need to prove the uniqueness. Assume
that $\boldsymbol{u}=0$ in \rf{4.1}$\sim$\rf{4.3}, we need to
prove $\pi_h\boldsymbol{u}=0$. From \rf{4.1}--\rf{4.2}, we have
$\pi_h\boldsymbol{u}\cdot\bo{n}=0$ on $\partial K$ and
\begin{eqnarray*}
(\hbox{div}\pi_h\boldsymbol{u},q)_K=-(\pi_h\boldsymbol{u},\nabla
q)_K+\int_{\partial K}\pi_h\boldsymbol{u}\cdot
\bo{n}qds=0,\;\forall\,q\in P_{k}(K).
\end{eqnarray*}
This implies $\hbox{div}\pi_h\boldsymbol{u}=0$ on $K$. So there
exists a function $\phi\in P_{k+2}(K)$ so that
curl$\,\phi=\pi_h\boldsymbol{u}$ (see \cite{Girault}). Since the
tangential derivative $\partial_\tau\phi=\hbox{curl}\,\phi\cdot\bo{n}=\pi_h\boldsymbol{u}\cdot
\bo{n}=0$ on $\partial K$, so $\phi=\phi_0=constant$ on $\partial K$.
Let $p=\phi-\phi_0$. Then, $p\in P^0_{k+2}(K)$ and
curl$\,p=\pi_h\boldsymbol{u}$. Taking $q=p$ in \rf{4.3}, we obtain
$\|\pi_h\boldsymbol{u}\|_{0,K}=0$ so that $\pi_h\boldsymbol{u}=0$.
Next, we prove conclusions \rf{4.4}$\sim$\rf{4.6}. Equation \rf{4.4} comes
directly from the Green's formula and \rf{4.1}--\rf{4.2}. From the
solution representation of linear system of equations
\rf{4.1}$\sim$\rf{4.3}, it is easy to see that on the reference
element $\hat{K}$,
\begin{equation}
\|\hat{\pi}_h\hat{\boldsymbol{u}}\|_{0,\hat{K}}\leq
\hat{C}(\|\hat{\boldsymbol{u}}\|_{0,\hat{K}}+\|\hat{\boldsymbol{u}}\|_{0,\partial
\hat{K}})\leq\hat{C}(\|\hat{\boldsymbol{u}}\|_{0,\hat{K}}+\|\nabla\hat{\boldsymbol{u}}\|_{0,\hat{K}}),\la{4.7}
\end{equation}
where we have used the trace inequality. Then, \rf{4.5} follows
from \rf{4.7} and a scale argument between $\hat{K}$ and $K$.
From \rf{4.5} and the unique existence, we also obtain
$$
\pi_h\bo{u}=\bo{u},\;\forall\,\bo{u}\in [P_{k+1}(K)]^2,\;\;\|\pi_h\boldsymbol{u}\|_{0,K}\leq
C\|\boldsymbol{u}\|_{1+s,K},\, 0\leq s\leq k+1.
$$
Hence, estimate \rf{4.6} can be derived by using the
Bramble-Hilbert lemma.\zb

Using the properties of operator $\pi_h$ and the argument of Lemma 3.2, we first can prove the discrete embedding inequality:
\begin{equation}
\|v^0\|\leq C\|\nabla_{w}v\|_h,\;\forall\, v\in
S_h^0,\la{4.8}
\end{equation}
and then prove the following result (see Theorem 3.1).\\
\dl{ 4.1}\quad{\em Let $V_h\times M_h$ be the weak Galerkin finite element space pair with $k\geq 1$. Then, the solution $(\bo{u}_h,p_h)$ of weak Galerkin finite element equations \rf{2.17}--\rf{2.18} uniquely exists and satisfies the stability estimate}
\begin{equation}
\|\nabla_{w}\bo{u}_h\|_h+\|p_h\|\leq
C\|\bo{f}\|.\la{4.9}
\end{equation}

By means of this extended projection $\pi_h\bo{u}$ and the parallel arguments to those in Section 3, we can give the following theorem.\\
\dl{ 4.2}\quad {\em For $k\geq 1$, let $(\bo{u},p)\in [H^{2+k}(\Omega)]^2\times H^{1+k}(\Omega)$ and $(\bo{u}_h,p_h)\in V_h\times M_h$ be the solutions of Stokes problem \rf{2.1}-\rf{2.3} and weak Galerkin finite element equations \rf{2.17}-\rf{2.18}, respectively. Then, we have}
\begin{eqnarray}
&&\|\nabla\bo{u}-\nabla_w\bo{u}_h\|_h+\|p-p_h\|\leq Ch^{k+1}(\|\bo{u}\|_{k+2}+\|p\|_{k+1}),\la{4.14}\\
&&\|Q_h^0\bo{u}-u_h^0\|\leq Ch^{k+2}(\|\bo{u}\|_{k+2}+\|p\|_{k+1})+Ch\|\bo{f}-Q_h^0\bo{f}\|.\la{4.14}
\end{eqnarray}
In particular, if $\bo{f}\in [H^{k+1}(\Omega)]^2$, then
\begin{equation}
\|Q_h^0\bo{u}-u_h^0\|\leq Ch^{k+2}(\|\bo{u}\|_{k+2}+\|p\|_{k+1}+\|\bo{f}\|_{k+1}).
\end{equation}
\section{Numerical experiment}
\setcounter{section}{5} \setcounter{equation}{0} In this section,
we examine the performance of the weak Galerkin finite element method described in Section 2. We apply this method to the following Stokes problem
\begin{eqnarray*}
-\Delta \bo{u} + \nabla p = \bo{f}, & \hspace{0.2in} \text{in~~} \Omega,  \\
\hbox{div}\bo{u} = 0, & \hspace{0.2in} \text{in~~} \Omega,  \\
\bo{u} = \bo{g}, & \hspace{0.2in} \text{on~~} \partial \Omega,
\end{eqnarray*}
where $\Omega = [0, 1]^2$, functions $\bo{f}$ and $\bo{g}$ are chosen such that the exact solution to this problem is
\begin{eqnarray*}
\bo{u}(x, y) &=& [xcos(y), cos(x)-sin(y)], \\
p(x, y) &=& x^3y - y^3 + 1/8.
\end{eqnarray*}

In the numerical experiments, we first partition $\Omega=[0,1]^2$ into a
regular triangle mesh $T_h$ with mesh size $h=1/N$. Then, the
refined meshes are obtained by using the edge bisection partition. Thus, we obtain a
mesh series $T_{h/2^j}, j=0,1,\dots$. We examine the computation error
for velocity and pressure approximations in the discrete
$H^1$-norm  and the $L_2$-norm, respectively. The numerical convergence rate $r$ is computed by using the formula
$r=\ln(e_h/e_{\frac{h}{2}})/\ln 2$, where $e_h$ is the computation
error. Numerical results are given in Table 1 for $k=0$ and Table 2 for $k=1$ with successively halved mesh size $h$.
We observe that the discrete solutions have a good approximation accuracy and the convergence rates are consistent with the theoretical prediction. Moreover, we also find from the computation error $\|Q_h^0\bo{u-u}_h^0\|$ that $\bo{u}_h^0$ is superclose to the projection $Q_h^0\bo{u}$ of solution $\bo{u}$. Numerical experiments verify the effectiveness of this weak Galerkin finite element method for Stokes problem.

\section{Conclusion}
\setcounter{section}{6} \setcounter{equation}{0}
We present a weak Galerkin finite element method for solving Stokes problem. Compared with those known works \cite{Xie,Lin3,Wang3,Zhai}, the main feature of our method is that the proposed weak Galerkin finite element scheme is stable without adding any stabilized term or penalty term and the velocity and pressure space pair satisfies the discrete inf-sup condition. Using the discrete inf-sup condition and a weak embedding inequality established here, we derive the optimal error estimates in the $H^1$- and $L_2$-norms for velocity and in the $L_2$-norm for pressure, respectively. Another important element in favor our method is that the proposed velocity and pressure space pair also can be used for Navier-Stokes problem. This is our ongoing work.

\begin{table}
\begin{center}
{\bf Tabel 1\quad}{\small History of convergence for discrete velocity and pressure for $k=0$\\[0.2cm]} \label{tab:1}

\begin{tabular}{cccc}
\hline
&$\|\nabla_{w}\bo{u}_h-\nabla \bo{u}\|_h$&\quad$\|p-p_h\|$&$\|Q_h^0\bo{u} - \bo{u}_h^0\|$\\
mesh $h$&error \quad\quad rate&\qquad error \qquad\quad\ rate&\qquad error \qquad\quad \ rate\\
\hline
1/10&2.8934e-02\qquad \ \ \ -   &\quad  2.9406e-02 \qquad \ \ -     & \quad 6.5665e-04\qquad \ \ \ - \\
1/20&1.4587e-02\quad 0.98805 &\qquad 1.4666e-02\qquad 1.0036  &   \qquad 1.6732e-04\qquad 1.9725 \\
1/40&7.3118e-03\quad 0.99642 &\qquad 7.3244e-03\qquad 1.0017  &   \qquad 4.2078e-05\qquad 1.9915\\
1/80&3.6586e-03\quad 0.99895 &\qquad 3.6605e-03\qquad 1.0007  &   \qquad 1.0538e-05\qquad 1.9974\\
1/160&1.8297e-03\quad0.99970 &\qquad 1.8300e-03\qquad 1.0002  &   \qquad 2.6360e-06\qquad 1.9992\\
1/320&9.1489e-04\quad0.99991 &\qquad 9.1493e-04\qquad 1.0001  &   \qquad 6.5911e-07\qquad 1.9998\\
\hline
\end{tabular}
\end{center}
\end{table}

\begin{table}
\begin{center}
{\bf Tabel 2\quad}{\small History of convergence for discrete velocity and pressure for $k=1$\\[0.2cm]} \label{tab:1}

\begin{tabular}{cccc}
\hline
&$\|\nabla_{w}\bo{u}_h-\nabla \bo{u}\|_h$&\quad$\|p-p_h\|$&$\|Q_h^0\bo{u} - \bo{u}_h^0\|$\\
mesh $h$&error \quad\quad rate&\qquad error \qquad\quad\ rate&\qquad error \qquad\quad \ rate\\
\hline
1/10&1.1746e-03\qquad \ \ \ -   &\quad  1.1186e-03 \qquad \ \ -     & \quad 1.0988e-05\qquad \ \ \ - \\
1/20&2.9579e-04\quad 1.9896 &\qquad 2.7978e-04\qquad 1.9994  &   \qquad 1.3842e-06\qquad 2.9887 \\
1/40&7.4183e-05\quad 1.9954 &\qquad 6.9969e-05\qquad 1.9995  &   \qquad 1.7377e-07\qquad 2.9938\\
1/80&1.8573e-05\quad 1.9979 &\qquad 1.7496e-05\qquad 1.9997  &   \qquad 2.1772e-08\qquad 2.9966\\
1/160&4.6466e-06\quad1.9990 &\qquad 4.3746e-06\qquad 1.9998  &   \qquad 2.7294e-09\qquad 2.9982\\
1/320&1.1621e-06\quad1.9995 &\qquad 1.0937e-06\qquad 1.9999  &   \qquad 3.4084e-10\qquad 2.9991\\
\hline
\end{tabular}
\end{center}
\end{table}
\vspace{0.2cm}

\section*{Acknowledgments}
The authors would like to thank the anonymous referees for many helpful suggestions which improved the
presentation of this paper. This work was supported by the National Natural Science Funds of
China, No. 11371081 and the State Key Laboratory of Synthetical
Automation for Process Industries Fundamental Research Funds, No.
2013ZCX02.

\end{document}